\input amstex
\input epsf
\input Amstex-document.sty
\define\map{\text{ map}}
\define\Diff{\text{ Diff}}
\define\hol{\text { hol}}
\define\lra{\longrightarrow}
\define\st{\text{ St}}
\define\th{\text{Th}}
\define\emb{\text{ Emb}}
\define\so{\Omega ^{\text{SO}}}

\pageno 447

\def\shorttitle{Strings and the Stable Cohomology of Mapping Class Groups}
\def\shortauthor{U. Tillmann}

\topmatter
\title\nofrills
{\boldHuge Strings and the Stable Cohomology of Mapping Class Groups}
\endtitle

\author \Large Ulrike Tillmann*
\endauthor

\thanks {*Math. Inst., 24--29 St. Giles, Oxford OX1 3LB, UK.
E-mail: tillmann\@maths.ox.ac.uk}
\endthanks

\abstract\nofrills 

\vskip 4.5mm


\endabstract

\endtopmatter
\document

\baselineskip 4.5mm
\parindent 8mm

\specialhead\noindent \boldLARGE
1. Introduction
\endspecialhead

Twenty years ago,  Mumford initiated the systematic study of the cohomology
ring of  moduli spaces of Riemann surfaces.
Around  the same time,
Harer proved that the homology of the mapping class groups
of oriented surfaces
is independent of the genus in low degrees, increasing with the genus.
The  (co)homology of mapping class groups thus stabelizes.
At least rationally, the mapping class groups have the same (co)homology
as the corresponding  moduli spaces.  This prompted Mumford  to
conjecture that the stable rational
cohomology of moduli spaces is generated by
certain tautological classes that  he defines. Much of the recent interest in
this subject is motivated by mathematical physics and, in particular,
by string theory.
The study of the  category of strings led to the discovery of
an infinite loop space, the cohomology of which is the
stable cohomology of the mapping class groups.
We  explain here  a homotopy theoretic approach to Mumford's conjecture
based on
this fact. As byproducts infinite families of torsion classes in the
stable cohomology are detected, and the divisibility of the   tautological
classes is determined.
An analysis of the category of strings in a background space leads to
the formulization of
a parametrized version of
Mumford's conjecture.

The paper is chiefly a summary of the author's work and her collaboration
with Ib Madsen. Earlier this year Madsen and Weiss announced a solution of
Mumford's conjecture. We  touch on some of the consequences and
the ideas behind this most
exciting new developement.

\specialhead\noindent\boldLARGE {2. Mumford's conjecture}
\endspecialhead

Let $F^s_{g,n}$ be an oriented, connected surface of genus $g$ with
$s$ marked points and $n$ boundary components. Let
$\Diff (F^s_{g,n})$ be  its group of orientaton preserving
diffeomorphisms
that fix the $n$ boundary components pointwise and permute the
$s$ marked points. By [2], for genus  at least 2,
$\Diff (F^s_{g,n})$ is homotopic to its group of
components, the mapping class group $\Gamma ^s_{g,n}$.
Furthermore,
if the surface has boundary,
$\Gamma ^s_{g,n}$ acts freely on Teichm\"uller space. Hence,
$$
B\Diff (F^s_{g,n}) \simeq B\Gamma ^s_{g,n} \simeq \Cal M ^s_{g,n}
\quad \text{ for } \, n\geq 1, g\geq 2,
$$
where $\Cal M^s_{g,n}$ denotes
the moduli space of Riemann surfaces appropriately
marked.
When $n=0$, the action of the mapping class group on Teichm\"uller space
has finite stabilizer groups and
the latter is only a  rational
equivalence.

We recall Harer's homology
stability theorem [4] which plays an important role through out
the paper.

\proclaim{Harer Stability Theorem 2.1 } $H_* B\Gamma ^s_{g,n}$ is independent of $g$ and $n$ in degrees $3* <
g-1$.
\endproclaim

Ivanov [5], [6] improved the stability range to $2* <g-1$ and proved a
version with twisted coefficients.
Glueing a torus with two boundary components to a surface $F_{g,1}$ induces a
homomorphism $ \Gamma _{g,1} \to \Gamma _{g+1,1}$. Let
$\Gamma _\infty := \lim _{g\to \infty} \Gamma _{g,1}$
be the stable mapping class
group.

Mumford [12]
introduced certain
tautological classes in the cohomology
of the moduli spaces $\Cal M _g$. Topological
analogues were studied by Miller [10] and Morita  [11]:
Let $E$ be the universal $F$-bundle over $B\Diff (F)$,  and let
$T^vE$ be its vertical
tangent bundle with
Euler class   $e \in H^2 E$.
Define
$$
\kappa_i:= \int _F e^{i+1} \quad \in \, H^{2i} B\Diff (F).
$$
Here $\int_F$ denotes \lq \lq integration over the fiber" - the Gysin map.
Miller and Morita
showed independently
that the rational cohomology of the stable mapping class group
contains the polynomial ring on the $\kappa_i$.

\proclaim{Mumford Conjecture 2.2} $H^* (B\Gamma _\infty; \Bbb Q) \simeq \Bbb Q [\kappa_1, \kappa _2, \dots].$
\endproclaim

\subheading {2.1. Remark}

The stable cohomology of the decorated mapping class groups is known modulo $H^* B\Gamma_\infty$ because of
decoupling  [1]. For example, let $\Gamma ^s_{\infty}:= \lim \Gamma ^s_{g,1}$. The following is a consequence of
Theorem 2.1.

\proclaim{Proposition 2.3} $(B\Gamma ^s_{\infty} )^+ \simeq B\Gamma ^+_\infty \times B(\Sigma _s \wr S^1) ^+.$
\endproclaim

Here $Y^+$ denotes Quillen's plus-construction on $Y$ with respect to the
maximal perfect subgroup of the fundamental group.
It is important to note that {\bf  the plus construction
does not change the (co)homology}. In particular,
$$
H^* B\Gamma _\infty = H^* B\Gamma ^+_\infty.
$$

\specialhead\noindent\boldLARGE {3. String category}
\endspecialhead

The category underlying the quantum mechanics of a state space $X$
is the path category $\Cal P X$.
Its objects are particles represented by
points in $X$. As time evolves a particle sweeps
out a path. Thus  a morphism between particles $a$ and $b$ is a
continuous path  in $X$
starting at $a$ and ending at $b$.
Concatenation of paths defines the composition in the category.
$$
\Cal P X = \left\{\aligned \text{objects }:\quad &  a, b, \dots \in X,   \\
\text{morphisms }: \quad & \coprod _{t>0} \map ([0,t], X). \endaligned \right.
$$

In string theory, the point
objects are replaced by closed loops in $X$. As time evolves these strings
sweep out a surface. Thus the space of morphisms from one string to
another is now described by a continuous
map from an oriented surface $F$ to $X$. The parametrization of the path
should be
immaterial. To reflect this,  take
homotopy orbits under the action
of $\Diff (F)$.
\footnote { Strings should also be independent of the parametrization. One
should therefore
take homotopy orbit spaces of the objects under the $S^1$ action.
In that case the diffeomorphisms of the surfaces need not be the identity
on the boundary. The resulting category has  the same homotopy
type as $\Cal S X$ in the sense that its classifying space is homotopic
to that of $\Cal SX$.}
Composition is given  by  concatenation of
paths, i.e. by glueing of surfaces along outgoing
and incoming boundary components.

To be more precise, let $LX=\map( S^1,X)$ denote the free loop space
on $X$.
A cobordism $F$
is a finite
union
$$
F_{g_1, n_1 +m_1} \cup \dots \cup F_{g_k, n_k + m_k}
$$
where $n=\Sigma _i n_i $ boundary components are considered
incoming and $m= \Sigma_i m_i$ outgoing.
For technical reasons we will assume $m_i >0$.
The category of strings in $X$ is then
$$
\Cal S X = \left\{\aligned \text{objects }: &\quad  \alpha, \beta,  \dots \in
\coprod_{n\geq 0}( LX)^n, \\
\text{morphisms }: &\quad \coprod _{F} E\Diff (F)\times _{\Diff (F)} \map (F, X). \endaligned \right.
$$
The disjoint union is taken over all cobordisms $F$, one for each
topological type.

\subheading{3.1. Elliptic elements}

The category $\Cal S X$ was first introduced by  Segal [14]. A functor from the path category $\Cal P X$ to the
category of $n$-dimensional vector spaces and their isomorphisms defines a vector bundle on $X$ with connection.
In particular, it defines an element in the $K$-theory of $X$. A functor from $\Cal S X$ to an appropriate
(infinite dimensional) vector space category is also referred to as a gerbe (or $B$-field) with connection. In
[14], Segal proposes this as the underlying geometric object of elliptic cohomology. More recently, this notion
has been refined by Teichner and Stolz.

\subheading {3.2. Conformal field theory}

The category $\Cal S := \Cal S (* )$ is studied in conformal field theory [15]. Its objects are the natural
numbers and its morphims are Riemann surfaces. A conformal field theory (CFT) is a linear space $\Bbb H$ which is
an algebra over $\Cal S $. Thus each element in $\Cal M_{g, n+m} $ defines a linear map from $\Bbb H ^{\otimes n}$
to $\Bbb H ^{\otimes m}$. The physical states of a  topological conformal field theory (TCFT) form a  graded
vector space $A_*$. Each element of the homology $H_* \Cal M _{g, n+m}$ defines a linear map from $A^{\otimes
n}_*$ to $A_*^{\otimes m}$.

\subheading {3.3. Gromov-Witten theory}

Let $X$ be a symplectic manifold. A model for the homotopy orbit spaces
in the definition of  $\Cal SX$ is the fiber bundle
$\map (F_{g,n}, X)   \to \Cal M_{g,n} (X)
\to \Cal M_{g,n}$
over the Riemann moduli space. In each fiber,
$F$ comes equipped with a complex structure and we may replace
the continuous maps by the space of pseudo-holomorphic maps
$\hol (F_{g}^n, X)$ yielding a category $\Cal S^{\hol} X$.
This is the category relevant to Gromov-Witten theory. Note,
for $X$ a complex Grassmannian, a generalized
flag manifold, or a loop group,  the degree
$d$-component of
$\hol (F_g, X)$ approximates
the  components of $\map (F_g, X)$ in homology.
The categories $\Cal SX$ and $\Cal S^{\hol} X$
are therefore closely related.

\specialhead\noindent \boldLARGE 4. From categories to infinite loop spaces
\endspecialhead

There is a functorial way to associate to a category $\Cal C$ a topological
space $|\Cal C|$, the realization of its nerve. It takes equivalences of
categories to
homotopy equivalences. It  is a generalization of the
classifying space construction of a group:
$|G| = BG$ where $G$ is identified with the category of a single object and
endomorphism set $G$.
The path-category $\Cal PX$ of a connected space $X$ is a many object group
up to homotopy. The underlying \lq \lq group" is the space
$\Omega X =\map _*(S^1, X)$ of based loops in $X$.
Again one has $|\Cal PX| \simeq B (\Omega X) \simeq X.$
A functor from $\Cal PX$ to $n$-dimensional vector spaces
and their isomorphisms thus defines a map
$$
X \lra BGL_n \Bbb C,
$$
and hence an element in the $K$-theory of $X$.
Motivated by this, we would like to understand the classifying space
of the string category $\Cal S X$ and its relation to elliptic cohomology.

\proclaim{Definition 4.1} $\st (X) := \Omega | \Cal S X|.$
\endproclaim

\proclaim{Theorem 4.2}
$\st$ is a homotopy functor from the category of topological spaces to
the category of infinite loop spaces.
\endproclaim

We recall that $Z$ is an infinite loop space if it is homotopic to some
$Z_0$ such that successive based spaces $Z_i$ can be found
with homeomorphisms $\gamma _i: Z_i \simeq \Omega Z_{i+1}$. Any infinite
loop space $Z$ gives rise to a generalized homology theory $h_*$ which
evaluated on a space $Y$ is
$$
h_* Y := \pi _* \lim _{i \to \infty} \Omega ^i (Z_i \wedge Y).
$$
Infinite loop spaces are abelian groups up to homotopy in the strongest sense.

The proof of Theorem 4.2 can be sketched as follows, compare [16].
$\Cal SX$ is a symmetric monoidal category under disjoint union.
Infinite loop space machinery (see for example [13])
implies that its classifying space $|\Cal SX|$ is
a homotopy
abelian  monoid in the strongest sense.
But $\pi_0 |\Cal SX|
= H_1 X$ is a group. Hence homotopy inverses exist and  $|\Cal SX|$ is
an infinite loop space.

Using another piece of infinite loop space machinery
(a generalization of the  group completion theorem)   and
Harer Stability Theorem 2.1,
one can identify the string theory of a point as
$ \Bbb Z \times B\Gamma ^+_\infty$.
As an immediate consequence
we have

\proclaim{ Corollary 4.3} ([16]) $\st (*) \simeq \Bbb Z \times B\Gamma^+_\infty$ is an infinite loop space.
\endproclaim

\specialhead \noindent \boldLARGE {5. CFT-operad}
\endspecialhead

We offer now a different perspective on
Theorem 4.1 and Corollary 4.3.
Let $\Cal M = \{ \Cal M_n\} _{n\geq 0}$ with $\Cal M(n) = \coprod _{g\geq 0}
B\Gamma _{g,n+1}$ be the operad contained in the CFT category $\Cal S$.
A  space  $X$ is an algebra over $\Cal M$ if there are compatible maps
$\Cal M (n) \times X^n \to X$. In particular $X$ has a monoid structure.
Let $\Cal G X$ be its  group completion. $\Cal GX$ is homotopic
to $X$ if and only if $\pi_0X$ is a group.

\proclaim {Theorem 5.1} ([18]) If $X$ is an algebra over $\Cal M$ then its group completion $\Cal GX$ is an
infinite loop space.
\endproclaim

CFT is therefore closely linked to the theory of infinite loop spaces.
The crucial point of the proof is a decoupling
result similar to Proposition 2.3.
The corresponding statement for TCFT's  implies that
Getzler's Batalin-Vilkovisky algebra structure on  the
physical states $A_*$ is stably trivial, see [17].
The following examples illustrate the strength of Theorem 5.1.

\subheading{Example 5.1}

Let $X_1$ be the disjoint union $\coprod _{g\geq 0} B\Gamma _{g,1}$. It has  a product induced by glueing
$F_{g,1}$ and $F_{h,1}$ to the \lq \lq legs" of a pairs of pants surface $F_{0,3}$. Miller observes in [10]  that
this induces a double loop space structure on the group completion $\Cal G X_1 \simeq  \Bbb Z \times B\Gamma
^+_\infty$. This product extends to an $\Cal M$-algebra structure. Hence Miller's double loop space structure
extends to an infinite loop space structure. Wahl [19] proved that it is equivalent to the infinite loop spaces
structure implied by Corollary 4.3.

\subheading{Example 5.2}

Let $X_2$ be the disjoint union of the Borel constructions
$$
E_g:= E\Diff (F_{g,1}) \times _{\Diff (F_{g,1})} \map (F_{g,1}, \partial ; X,*).
$$
As the functions to $X$ map the boundary to a point, they can be extended
from $F_{g,1}$ to
$F_{g+1,1}$ via the constant map.
$X_2$ thus becomes
an $\Cal M$-algebra and  $\Cal G X_2 = \Bbb Z \times (\lim _{g\to
\infty} E_g)^+$ is an infinite loop space.
$\Cal GX_2$ is homotopic to $\st (X)$ when $X$ is simply connected.

\subheading{Example 5.3}

Let $C(F_{g,1};X)$ denote the space of unordered
configurations  in the interior of $F_{g,1}$ with labels in
a connected space $X$. Let   $C_g$ be its Borel construction.
Their disjoint union  defines an $\Cal M$-algebra $X_3$.
The following  decoupling result determines its group completion.

\proclaim {Proposition 5.2} $\Cal G X_3 \simeq \Bbb Z \times (\lim _{g \to \infty} C_g )^+ \simeq \Bbb Z \times
B\Gamma ^+_\infty \times Q (BS^1 \wedge  X_+) .$
\endproclaim

$Q =\lim \Omega^\infty S^\infty$ is the free infinite loop space
functor and $X_+$ denotes $X$ with a disjoint basepoint.
Note the close relation with the above example.
By work of McDuff and B\"odigheimer, there is a homotopy equivalence
$$
C(F_{g,1}; X) \simeq \map (F_{g,1}, \partial; S^2 \wedge X).
$$
Note though that
the induced $\Diff (F_{g,1})$-action from the left on the right is non-trivial
on the sphere in the target space.

\specialhead \noindent \boldLARGE {6. Refinement of Mumford's conjecture}
\endspecialhead

Infinite loop spaces are relatively rare and the
question arises whether $\Bbb Z\times B\Gamma ^+_\infty$ can be
understood in terms of  well-known infinite loop spaces.
This question was addressed in joint work with Madsen.

Let $\Bbb P^l$ be the Grassmannian of oriented
2-planes in $\Bbb R^{l+2}$ and let $-L_l$ be the complement of the
canonical 2-plane bundle $L$ over $\Bbb P^l$. The one-point compactification,
the Thom space $\th (-L_l)$, restricts on the subspace $\Bbb P^{l-1}$ to
the suspension of $\th (-L_{l-1})$. Taking adjoints yields maps
$\th (-L_{l-1}) \to \Omega \th (-L_l)$, and we may define
$$
\Omega ^\infty \th (-L) := \lim _{l\to \infty} \Omega ^l \th (-L_l).
$$
More generally,  for any space $X$,
define
$$
\Omega ^\infty (\th (-L) \wedge X_+) := \lim _{l\to \infty} \Omega ^l
(\th (-L_l) \wedge X_+).
$$

\proclaim {Conjecture 6.1} There is a homotopy equivalence of infinite loop space
$$
\alpha : \st (X) \overset \simeq \to \lra
\Omega ^\infty (\th (-L) \wedge X_+).
$$
\endproclaim

\subheading{Remarks 6.2}

For $X=*$, Conjecture 6.1 postulates a homotopy equivalence $ \alpha: \Bbb Z\times B\Gamma ^+_\infty \to \Omega
^\infty \th (-L). $ A proof of this has been announced by Madsen and Weiss, see Section 8. The Mumford conjecture
follows from this as we will explain presently. Conjecture 6.1  claims in addition that  $\st (\_ ) $ is a
homology functor, i.e. $\pi_* \st (\_ )$ is the homology theory associated to the infinite loop space $\Omega
^\infty \th (-L)$.

The infinite loop space $\Omega ^\infty \th (-L)$ is well-studied, more recently because of its relation to
Waldhausen $K$-theory. The inclusion of $-L_l$ into the trivial bundle $  L\oplus (-L_l) \simeq \Bbb P^l \times
\Bbb R^{l+2}$ induces a map
$$
\omega: \Omega ^\infty \th (-L) \lra Q (\Bbb P^\infty _+).
$$
$\omega $ has homotopy fibre $\Omega ^2 Q (S^0)$. As the stable homotopy groups
of the sphere are torsion in positive dimensions, $\omega$ is a rational
equivalence.
Let $\Bbb P^\infty \to BU$ be the map that classifies $L$.
By Bott periodicity, the map can
be extended to the free infinite loop space $L: Q_0 (\Bbb P^\infty) \to
BU$. The subscript here indicates the 0-component.
$L$ has a splitting and is well known to be a rational
equivalence:
$$
H^* (\Omega_0 ^\infty \th (-L); \Bbb Q) \overset {\omega ^* \circ L^*} \to
\simeq \Bbb Q [ c_1, c_2, \dots ].
$$
The $\Bbb Z/p$-homology of $\Omega ^\infty \th (-L)$
has recently been determined by Galatius [3].

\vskip .2in \centerline{\hbox{\hskip 0.2cm \epsfysize=3cm\epsfbox{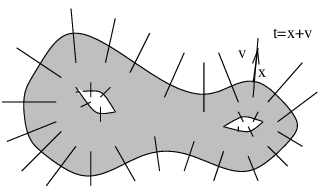}}}
\hfil {Figure 1: Surface $h(F)\subset \Bbb R ^{l+2}$ with tubular neighborhood
$U$.}
\hfil

\vskip .1in
To define a map $\alpha $ comes down to defining
maps from the morphism spaces of
$\Cal S X$. (See also Example 5.2.)
$\alpha$  is the homotopy theoretic
interpretation of the formula defining  $\kappa_i$
where the wrong way map $\int _F$ is
replaced by the (pre)transfer map of Becker and Gottlieb.
We  give now an explicid description of this map.

For simplicity, let $F$ be a closed surface. Consider the space of
smooth embeddings $\emb (F, \Bbb R^{l+2})$. By Whitney's embedding theorem,
as $l \to \infty$,
it may serve as a
model for $E\Diff (F)$.
Let
$$
(h, f) \in \emb (F, \Bbb R ^{l+2}) \times _{\Diff(F)} \map (F,X).
$$
Choose a tubular neighborhood $U$ of $h(F)$ such that
every  $t \in U$
can uniquely be written as $x+v$ with $x \in h(F)$ and $v$ normal to the
tangent plane $T_xh(F)$. $\alpha$ sends $(h,f)$ to the continuous function
$\alpha (h,f)
: S^{l+2} \to \th (-L_l) \wedge X_+$   defined by
$$
t \quad  \mapsto \quad  \left\{ \aligned \infty \quad \quad  &\text{ if }
t\notin U, \\
((T_xh(F), v), f(h^{-1} (x)) \quad  &\text{ if } t=x+v \in U. \endaligned \right.
$$

In [8], $\alpha$ is shown to be a 3-connected map of infinite loop spaces and
the tautological classes are identified.
Let
$i! (ch_i) \in H^{2i} BU$ denote the $i$-th integral Chern character class.
Then
$$
\kappa_i = \alpha ^*\circ \omega ^* \circ L^* (i! ch_i).
$$

\specialhead \noindent \boldLARGE {7. Splittings and (co)homological results}
\endspecialhead

The main result of [8] is a partial splitting of the composition
$$
\omega \circ \alpha :
\Bbb Z \times B\Gamma ^+_\infty \lra Q (\Bbb P^\infty _+) \simeq Q (S^0) \times
Q (\Bbb P ^\infty).
$$
This is achieved by constructing a map $\mu$ from $\Bbb P _+$ to $\Bbb Z \times
B\Gamma ^+_\infty$ and then extending it to the free infinite loop space
$Q (\Bbb P^\infty_+)$ utilizing the infinite loop space
structure on $\Bbb Z\times B\Gamma ^+_\infty$. In order to construct $\mu$,
approximate $\Bbb P ^\infty \simeq BS^1$  by the classifying spaces
of cyclic groups $C_{p^n}$ for $n\to \infty$,
one prime $p$ at a time, as the cyclic groups can
be mapped into suitable mapping class groups.
However, this means that we have
to work with  $p$-completions.

Let $Y^\wedge _p$ denote the $p$-completion of $Y$ and
$g \in \Bbb Z ^\times _p$ be a topological generator of the $p$-adic units
($g=3$ if $p=2$).
Denote by $\psi ^k: \Bbb P^\infty \to (\Bbb P ^\infty) ^\wedge _p $
the map that represents
$k$ times the first Chern class in $H^2 (\Bbb P ^\infty , \Bbb Z _p)$.

\proclaim {Theorem 7.1. [8]} There exists a map $\mu: (Q(S^0 ) \times Q (\Bbb P^\infty ))^\wedge _p \to (\Bbb Z
\times B\Gamma ^+_\infty)^\wedge _p$ such that
$$
\omega \circ \alpha \circ \mu \simeq \left( \matrix
-2 & *  \\
0  & 1- g \psi^g
\endmatrix \right).
$$
\endproclaim

The map $1-g\psi^g$ induces multiplication by $1-g^{n+1}$ on $H_{2n} (\Bbb P^
\infty ; \Bbb Z_p)$ which is a $p$-adic unit precisely if
$n \neq -1 (\mod p-1)$.
The following applications of Theorem 7.1 are also found in [8].
There is a splitting  $Q(\Bbb P^\infty _+)^\wedge _p \simeq E_0 \times \dots
\times E_{p-2}$
corresponding to the idempotent decomposition of
$\Bbb Z[ \Bbb Z /p ^\times] \subset \Bbb Z_p [ \Bbb Z_p^\times]$.

\proclaim {Corollary 7.2} For some $W_p$, there is a splitting of infinite loop spaces
$$ (\Bbb Z \times B\Gamma ^+_\infty )^\wedge_p \simeq E_0 \times \dots \times
E_{p-3} \times W_p.
$$
\endproclaim

The $\Bbb Z/ p$-homology of $Q (\Bbb P^\infty_+)$ is well-understood in
terms of Dyer-Lashof operation. These are homology operations for infinite
loop spaces that are formally similar to the Steenrod operations.  For
each  generator $a_i \in H_{2i} \Bbb P^\infty =\Bbb Z $
there is an infinite family of $\Bbb Z/p$-homology classes
freely generated by the Dyer-Lashof operations.
The product $E_0 \times \dots \times E_{p-3}$
contains precisely those families for which $i\neq -1
(\text {mod } p-1)$, giving a huge collection of new $p$-torsion in
$H_* B\Gamma _\infty$.

For odd primes $p$, Madsen and Schlichtkrull [MS] found  split surjective maps
$l_0$ and $l_{-1}$ of infinite loop spaces such that the following
diagram is commutative
$$
\CD
\Omega ^\infty \th (-L) ^\wedge _p  @> \omega >>    Q (\Bbb P ^\infty _+ )
    ^\wedge _p  \\
@V l_{-1} VV    @V l_0 VV   \\
(\Bbb Z \times BU ) ^\wedge _p @> 1- g \psi ^g >>   (\Bbb Z \times BU )
    ^\wedge _p.
\endCD
$$

\proclaim{Corollary 7.3} For odd primes $p$ and some space $V_p$, there is a splitting of spaces
$$
(B\Gamma ^+_\infty ) ^\wedge _p \simeq  BU ^\wedge _p \times V_p.
$$
\endproclaim

This gives  a $\Bbb Z_p$-integral version of Miller and Morita's theorem:
the polynomial algebra $\Bbb Z_p [ c_1, c_2, \dots ]$ is a split summand
of $H^* (B\Gamma _\infty ; \Bbb Z _p)$.
The divisibility of the tautological classes $\kappa _i$
at odd primes $p$ can also be deduced from the above diagram.

\proclaim{Corollary 7.4} If $i = -1 (\text {mod } p-1)$, then $\kappa _i$ is divisible by $p ^{1+\nu_p (i+1)}$
where $\nu_p$ is the $p$-adic valuation. Otherwise, $p$ does not divide $\kappa _i$.
\endproclaim
In the light of [9] this result is sharp.

\specialhead \noindent \boldLARGE { 8. Geometric interpretation}
\endspecialhead

$\alpha: B\Gamma ^+_\infty \to \Omega ^\infty _0 \th (-L)$
is a homotopy equivalence if and only if it  induces an
isomorphism in oriented cobordism theory $\so_* $.
An element in $\so_n(B \Gamma ^+_\infty)
= \so_n (B\Gamma _\infty)$ is a cobordism class of
oriented surface bundles  $F \to E^{n+2} \overset \pi \to
\to M ^n$.
An element in $\so _n (\Omega ^\infty _0 \th (-L))$ is
a cobordism class of pairs $[ \pi: E ^{n+2} \to M^n, \hat \pi
]$ of smooth maps $\pi$ and stable bundle surjections from $TE$ to $\pi ^*TM$.
(Upto cobordism one can assume that $\hat \pi $ is a vector bundle
surjection.) $\alpha$ maps a bundle $[F\to  E\overset \pi \to \to M]$
to the pair
$[\pi: E \to M, D\pi]$ where $D\pi$ denotes the differential of $\pi$.
Hence, $\alpha$ is a homotopy equivalence if and only if
each cobordism class of pairs $[ \pi :E ^{n+2} \to M^n, \hat \pi]$
contains a \lq\lq unique" representative with $\pi$ a submersion.

It is this geometric formulation that underpins the solution to the Mumford
conjecture by Madsen and Weiss.
A key ingredient of the proof is the Phillips-Gromov
$h$-principle of submersion
theory: A pair $(g: X\to M, \hat g: TX \to g^*TM)$
can be deformed to a submersion
-- provided $X$ is open.
$E$ above, however, is closed.
The approach taken in [9] is to replace
 $\pi: E\to M$ by $g = \pi \circ pr_1
: X= E\times \Bbb R \to M$. Now the submersion
$h$-principle applies and $g$ can be replaced
by a submersion $f$. The proof then consists of a careful
analysis of the
singularities of the projection $pr_1: X\to \Bbb R$ on the fibers of
$f$.
At a critical point it uses
Harer's Stability
Theorem 2.1.

\proclaim {Madsen-Weiss Theorem 8.1} The map $\alpha : \Bbb Z \times B\Gamma ^+_\infty \to \Omega ^\infty \th
(-L)$ is a homotopy equivalence.
\endproclaim

\widestnumber\key{AA}

 \specialhead \noindent \boldLARGE {References}
\endspecialhead

\ref
\key 1\by C.-F. B\"odigheimer \& U. Tillmann \paper
{ Stripping and splitting
decorated mapping class groups}
\jour  Birkh\"auser, Progress in Math. \vol { 196}
\yr 2001 \pages 47-- 57
\endref

\ref
\key 2 \by C.J.  Earle \& J. Eells \paper { A fibre bundle description of
Teichm\"uller theory} \jour J. Diff. Geom. \vol 3 \yr 1969\pages 19--43
\endref

\ref
\key 3\by S. Galatius \paper Homology of $\Omega ^\infty \Sigma
\Bbb CP^\infty _{-1} $ and $ \Omega ^\infty \Bbb CP ^\infty _{-1}$
\jour preprint 2002
\endref

\ref
\key 4 \by J.L. Harer \paper
 Stability of the homology of the mapping class groups of
orientable surfaces \jour  Annals Math. \vol 121 \yr 1985
\pages 215--249
\endref

\ref
\key 5 \by N.V. Ivanov \paper
Stabilization of the homology of Teichm\"uller modular
groups
\jour
Lenin\-grad Math. J. \vol 1 \yr 1990 \pages 675--691
\endref

\ref
\key 6 \by N.V. Ivanov \paper
On the homology stability  for Teichm\"uller mudular
groups: closed surfaces and twisted coefficients
\jour  Mapping Class Groups and Moduli Spaces of Riemann Surfaces,
Contemp. Math. \vol 150 \yr 1993 \pages 149--194
\endref

\ref
\key 7 \by I. Madsen \& C. Schlichtkrull \paper
The circle transfer and $K$-theory \jour
AMS Contemporary Math.  \vol 258 \yr 2000 \pages 307--328
\endref

\ref
\key
8\by I. Madsen \& U. Tillmann \paper The stable mapping class group and
$Q (\Bbb CP ^\infty )$ \jour Invent. Math. \vol 145 \yr 2001 \pages 509--544
\endref

\ref
\key 9 \by I. Madsen \& M. Weiss \paper
Cohomology of the stable mapping class group
\jour in preparation
\endref

\ref
\key 10 \by E.Y. Miller \paper
The homology of the mapping class group \jour J. Diff. Geom.
\vol 24 \yr 1986 \pages 1--14
\endref

\ref
\key
11 \by S. Morita \paper
Characteristic classes of surface bundles \jour Invent. Math.
\vol 90 \yr 1987 \pages 551-577
\endref

\ref
\key 12 \by D. Mumford \paper
Towards an enumerative geometry of the moduli space of
curves \jour  Arithmetic  and Geometry, M. Artin and J. Tate, editors,
Progr.
Math., Birkhauser  \vol 36  \yr 1983 \pages 271--328
\endref

\ref
\key 13 \by G. Segal \paper Categories and cohomology theories
\jour Topology \vol 13
\yr 1974 \pages 293--312
\endref

\ref
\key 14 \by G. Segal
\paper Elliptic cohomology (after Landweber-Stong, Ochanine, Witten,
and others) \jour Seminar Bourbaki,  Asterisque \vol 161-162
\yr 1989 \pages 187--201
\endref

\ref
\key 15 \by G. Segal \paper
The definition of conformal field theory \jour manuscript
\endref

\ref \key
16 \by U. Tillmann \paper On the homotopy of the stable mapping class group
\jour Invent. Math. \vol 130 \yr 1997 \pages 257--275
\endref

\ref
\key 17 \by U. Tillmann \paper
Vanishing of the Batalin-Vilkovisky algebra structure
for TCFTs \jour Commun. Math. Phys \vol  205 \yr 1999 \pages 283--286
\endref

\ref
\key 18 \by U. Tillmann \paper
Higher genus surface operad detects infinite loop
spaces \jour Math. Ann. \vol  317 \yr 2000\pages 613--628
\endref

\ref
\key 19 \by N. Wahl \paper  Infinite loop space structure(s) on the stable
mapping class group \jour Oxford Thesis 2001
\endref

\enddocument